# Electromagnetic Pulse Propagation over Nonuniform Earth Surface: Numerical Simulation

**Alexei V. Popov and Vladimir V. Kopeikin**

Pushkov Institute of Terrestrial Magnetism,
Ionosphere and Radio Wave Propagation
*IZMIRAN, Troitsk, Moscow region, 142190 Russia*
*popov@izmiran.ru, kopeikin@izmiran.ru*

**1. Introduction**

Parabolic equation method proposed by Leontowich and Fock [1,2] is an efficient simulation approach to VHF propagation over the earth surface. Deep physical analysis and advanced mathematical methods [3,4] turned Leontovich's PE into a universal tool of diffraction theory. Its applications go far beyond the initial problem circle – e. g. [5-8]. The key role in this development played the decisive turn to straightforward numerical techniques pioneered by Malyuzhinets and Tappert [9,10].

In radio wave propagation, PE was used first to derive explicit analytical formulae for the EM field strength in model environments. A simplification has been reached by introducing the impedance boundary condition (BC) [11]. Taking into account tropospheric refraction ducts required the use of sophisticated asymptotic methods [12]. Further development (almost exclusively towards numerical implementation) was aimed at refined PE modifications [13-15], account for irregular terrain [16], introducing artificial transparent boundaries [17,18] and nonlocal BC to describe rough interfaces [19]. A non-stationary PE counterpart and a finite-difference (FD) scheme for its solution have been proposed by Claerbout and applied to seismic problems [13]. Afterwards, this "time-domain parabolic equation" (TDPE) was used to calculate acoustic propagation in ocean [20]. At the same time, little attempts of using TDPE to simulate EM pulse propagation in realistic environments are known.

In this paper, we consider computational aspects of EM pulse propagation along the nonuniform earth surface. For ultrawide-band pulses without carrier, TDPE results directly from the exact wave equation written in a narrow vicinity of the wave front. To solve it by finite differences, we introduce a time-domain analog of the impedance BC and a nonlocal BC of transparency reducing the open computational domain to a strip of finite width. Numerical examples demonstrate the influence of soil conductivity on the received pulse waveform which can be used in remote sensing.

For a high-frequency modulated EM pulse, TDPE arises as a convolution of PE solutions with the pulse envelope spectrum. In order to overcome computational difficulties, we develop an asymptotic approach based on the ray structure of the monochromatic wave field calculated at the carrier frequency. To accommodate complex-valued asymptotic solutions to the real initial condition we use the "analytic signal" approach introduced by Vainstein, Heyman and Felsen [21, 22]. An explicit solution of the time-domain transport equation reduces the computational procedure to numerical integration of standard PE at the carrier frequency and evaluation of a given 1D function in time domain. This diminishes computational expenses by 2-3 orders of magnitude and allows for pulsed wave field calculation in vast domains measured by tens of thousands wavelengths. As an example, we consider a problem of target altitude determination from overland radar data [23].

This work has been done in collaboration with the Institute for high-frequency technique (IHF), Stuttgart University. Preliminary results appeared as short papers [24,25], a Russian version has been published in [26]. We dedicate this publication to the memory of Leopold Benno Felsen.

**2. Monochromatic wave propagation**

Omitting technical details and method refinements – see [12, 19], recall PE based scheme of monochromatic wave propagation over a smoothly nonuniform earth surface $z = h(x)$ - Fig.1.

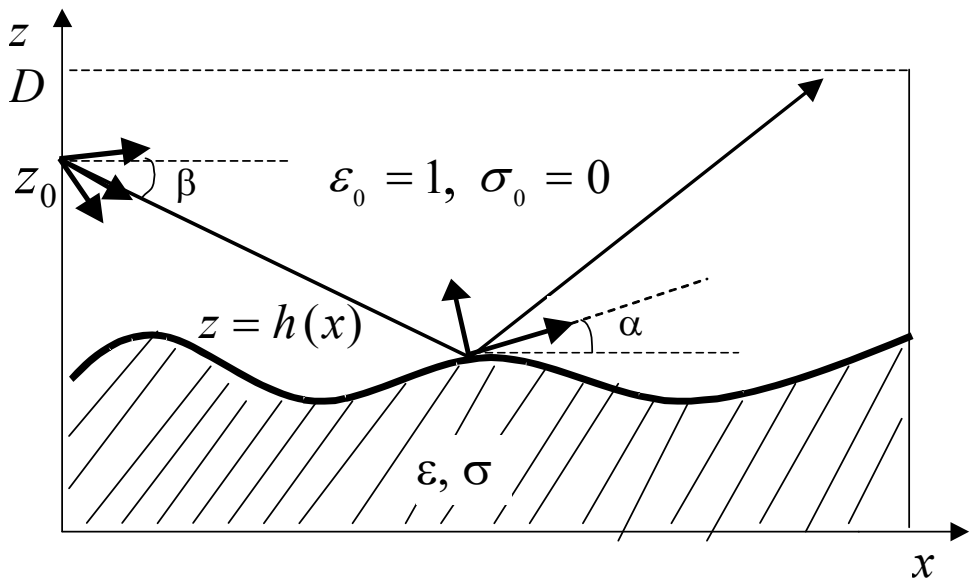

Fig. 1. Elevated source illuminating smoothly rolling terrain (sketch).

Horizontal magnetic component $H_y = H(x,z)$ satisfies Helmholtz equation

$$\frac{\partial^2 H}{\partial x^2} + \frac{\partial^2 H}{\partial z^2} + k^2 \tilde{\varepsilon} H = 0 \qquad (1)$$

with complex permittivity $\tilde{\varepsilon} = \varepsilon + 4\pi i \sigma/\omega$, where $\sigma$ is soil conductivity in the Gaussian units set. In the upper medium $\varepsilon = 1$, $\sigma = 0$ and the contact conditions at $z = h(x)$ are

$$H^+ = H^-, \qquad \frac{\partial H^+}{\partial n} = \frac{1}{\tilde{\varepsilon}} \frac{\partial H^-}{\partial n}. \qquad (2)$$

where $\partial/\partial n = \cos\alpha\, \partial/\partial z + \sin\alpha\, \partial/\partial x$, $\alpha = \arctan h'(x)$. At large distances from the source the wave field is sought as a plane wave with slowly varying complex amplitude:

$$H(x,z,t) \approx u(x,z,k) \exp[i(kx - \omega t)] \qquad (3)$$

Here, $k = \omega/c \equiv 2\pi/\lambda$ is the wave number, and the complex "attenuation function" $u(x,z,k)$ satisfies the Leontovich PE

$$2ik\frac{\partial u}{\partial x} + \frac{\partial^2 u}{\partial z^2} = 0, \qquad z > h(x) \qquad (4)$$

In 3D, divergence factor $1/\sqrt{x}$ must be added in (3). In this paper we use Gaussian initial condition

$$u(0,z,k) = \exp\left[\left(\frac{ik}{2\rho_0} - \frac{1}{w_0^2}\right)(z - z_0)^2 - ik\beta(z - z_0)\right] \qquad (5)$$

corresponding to an exact solution of PE (4)

$$u_i(x,z,k) = \sqrt{\frac{x_0}{x + x_0}} \exp\left\{ik\left[\frac{(z - z_0 - \beta x)^2}{2(x + x_0)} + \beta(z - z_0) - \frac{\beta^2}{2}x\right]\right\} \qquad (6)$$

- skewly propagating Gaussian beam with initial width $w_0$ and wave front radius $\rho_0$, determined by complex parameter $x_0 = \left(1/\rho_0 + 2i/k w_0^2\right)^{-1}$; $\beta$ being a small elevation angle.

Impedance approximation is based on wave beam contraction when entering a denser dielectric medium. Standard Leontovich BC [11]

$$\frac{\partial H^+}{\partial z} = -\frac{ik}{\sqrt{\tilde{\varepsilon}}} H^+, \qquad z = 0 \qquad (7)$$

results from the contact conditions (2) under the assumption of almost vertical propagation in the lower medium: $H^-(x,z) \approx T \exp(-ikz\sqrt{\tilde{\varepsilon}})$. For grazing angles (Fig. 2) this assumption breaks and a plane incident wave $H^+(x,z) = \exp[ik(x\cos\beta - z\cos\beta)]$ with small $|\beta| \ll 1$

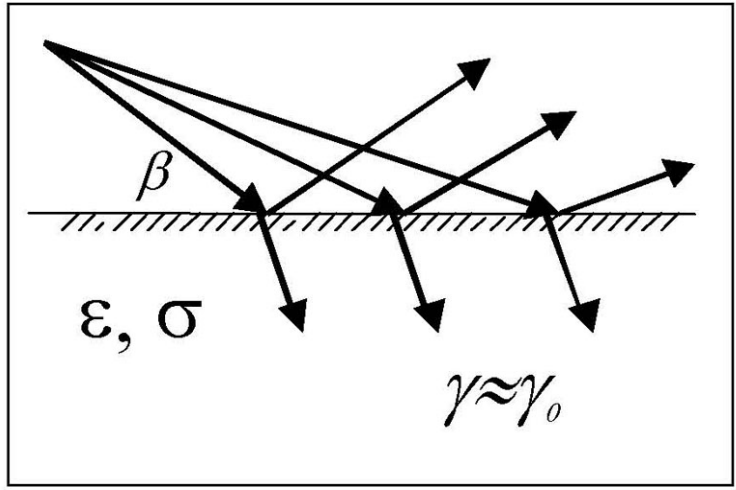

Fig. 2. Derivation of impedance BC for grazing angles.

enters the half-space $z < 0$ close to the total internal reflection angle $\gamma_0 = \arccos(\tilde{\varepsilon}^{1/2})$:

$\gamma \approx \gamma_0 + \dfrac{\beta^2}{2\sqrt{\tilde{\varepsilon}} \sin \gamma_0}$. Hence ensue

$$\begin{aligned} H^-(x,t) &= T(\gamma) \exp\left[ ik\sqrt{\tilde{\varepsilon}}(x \cos \gamma - z \sin \gamma) \right] \approx \\ &\approx T(\gamma_0) \exp\left[ ik(x - z\, tg\gamma_0) \right] \exp\left[ -\dfrac{i}{2} k\beta^2 (x + z\, ctg\gamma_0) \right] \equiv \\ &\equiv T(\gamma_0) \exp\left[ ik(x - z\, tg\gamma_0) \right] V(x + z\, tg\gamma_0) \end{aligned} \qquad (8)$$

In virtue of the superposition principle, Eq. (8) holds for an arbitrary paraxial wave packet with the corresponding slowly varying function $V(x)$. Eliminating the latter by differentiation and making use of (2), we obtain

$$\sqrt{\tilde{\varepsilon} - 1}\, \dfrac{\partial H^+}{\partial z} - \dfrac{1}{\tilde{\varepsilon}} \dfrac{\partial H^+}{\partial x} + ik\tilde{\varepsilon} H^+ = 0, \qquad z = 0 \qquad (9)$$

This modified impedance BC provides a more accurate approximation of the reflection coefficient, especially, in a vicinity of the Brewster angle $\beta_0 = \arcsin\left(\sqrt{\tilde{\varepsilon} + 1}\right)^{-1/2}$.

Fig. 3 allows one to compare the exact Fresnel reflection coefficient

$$R_F(\beta) = \dfrac{\tilde{\varepsilon} \sin \beta - \sqrt{\tilde{\varepsilon} - \cos^2 \beta}}{\tilde{\varepsilon} \sin \beta + \sqrt{\tilde{\varepsilon} - \cos^2 \beta}} \qquad (10a)$$

with the Leontovich approximation

$$R_L(\beta) = \dfrac{\sqrt{\tilde{\varepsilon}} \sin \beta - 1}{\sqrt{\tilde{\varepsilon}} \sin \beta + 1} \qquad (10b)$$

and that resulting from the modified impedance BC (9)

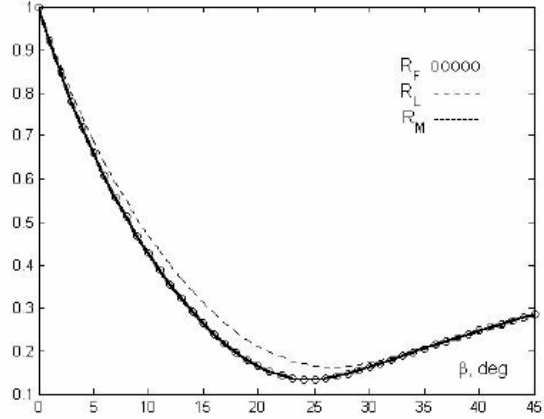

Fig.3. Comparison between exact Fresnel reflection coefficient and impedance approximations.

$$R_M(\beta) = \frac{\tilde{\varepsilon}\sin\beta - (\tilde{\varepsilon} - \cos\beta)/\sqrt{\tilde{\varepsilon} - 1}}{\tilde{\varepsilon}\sin\beta + (\tilde{\varepsilon} - \cos\beta)/\sqrt{\tilde{\varepsilon} - 1}} \qquad (10c)$$

Taking into account the boundary tilts $h'(x)$ and using "parabolic" approximation (3), we derive a modified BC for the attenuation function $u(x,z)$:

$$\frac{\partial u}{\partial z} + ik\left[\frac{\sqrt{\tilde{\varepsilon}-1}}{\tilde{\varepsilon}} - h'(x)\right]u = 0, \quad z = h(x) \qquad (11)$$

Contrary to the standard Leontovich BC (9), here it is not necessary to assume $|\tilde{\varepsilon}| \gg 1$ - Eq. (11) breaks only for $|\tilde{\varepsilon} - 1| \ll 1$ when nonlocal effects of wave interaction in both media are to be taken into account [19]. The impedance BC grants uniqueness of the solution of PE (4). In fact, calculating the energy flow through a vertical cross section one obtains

$$I(x) = \int_{h(x)}^{\infty} |u(x,z)|^2 dz, \qquad \frac{dI}{dx} = -\mathrm{Re}\left[\frac{\sqrt{\tilde{\varepsilon}-1}}{\tilde{\varepsilon}}\right] \cdot |u(x,h)|^2 \le 0 \qquad (12)$$

which proves stability and uniqueness of the boundary value problem solution.

Finite-difference methods of PE solution have been studied in early works by Malyuzhinets and coauthors [9, 27]. Further method development is described in monographs [5, 19]. We employ a six-point implicit FD scheme supplemented with the impedance BC (11) at $z = h(z)$ and a discrete approximation of the nonlocal transparency BC [17,28] imposed at the artificial computational boundary $z = z_{max}$:

$$\frac{\partial u}{\partial z}(x, z_{max}) = -\sqrt{\frac{2ik}{\pi}} \int_0^x \frac{\partial u}{\partial x}(\xi, z_{max}) \frac{d\xi}{\sqrt{x-\xi}} \qquad (13)$$

An example of simulated VHF propagation over irregular terrain is illustrated by Fig. 4.

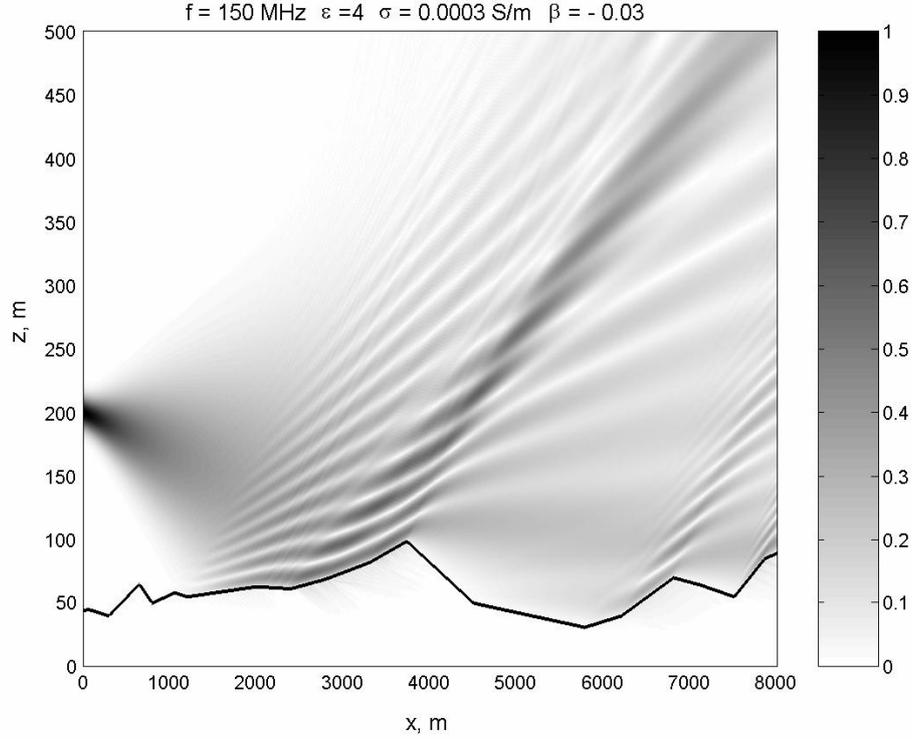

Fig. 4. VHF attenuation function over irregular earth surface.

## 3. Radio pulse propagation: Fourier synthesis

A straightforward way to describe EM transients is to convolve monochromatic wave fields with the signal spectrum. In a 1D case, the propagating pulse is a superposition of plane waves

$$H(x,t) = \frac{1}{2\pi} \int_{-\infty}^{\infty} \widetilde{F}(k) e^{i[kx - \omega(k)t]} dk \qquad (14)$$

In free space, $\omega = kc$, and formula (14) yields a dispersion-less traveling wave

$$H(x,t) = F(ct-x), \quad F(s) = \frac{1}{2\pi} \int_{-\infty}^{\infty} \widetilde{F}(k) e^{-iks} dk \qquad (15)$$

In a 2D environment, a natural generalization of the 1D solution (14) is a paraxial wave packet

$$H(x,z,t) = \frac{1}{2\pi} \int_{-\infty}^{\infty} \widetilde{F}(k) u(x,z,k) e^{ik(x-ct)} dk \qquad (16)$$

where $u(x,z,k)$ is a solution of the PE (4) at a fixed frequency $\omega = kc$.

The superposition (16) will approximate an exact solution of the wave equation if the spectrum $\widetilde{F}(x)$ is confined near a certain positive $k_0$ satisfying the PE applicability conditions: $k_0 \gg 2\pi/D \gg 2\pi/L$ where $D$ and $L$ are lateral and longitudinal characteristic scales of the problem. Consider a quasi-monochromatic pulse $F(ct) = f(ct)\cos\omega_0 t$ with duration $T \gg 2\pi/\omega_0$. Its Fourier transform

$$\widetilde{F}(k) = \frac{1}{2}\left[\widetilde{f}(k-k_0) + \widetilde{f}(k+k_0)\right], \tag{17}$$

where $\widetilde{f}(k) = c\int_{-\infty}^{\infty} f(ct)e^{i\omega t}dt$ is the envelope $f(ct)$ spectrum, contains negative frequencies not described by PE (4). Introducing complex signal $F_c(ct) = f(ct)\exp(-i\omega_0 t)$ eliminates the second term in (17). Still, the remaining "positive" component $\widetilde{f}(k-k_0)$ centered at $k_0$ may spread onto negative semi-axis. In order to avoid nonphysical effects of negative frequencies propagation, the "analytic signal" [21] can be used, defined as a one-side inverse Fourier transform of the truncated spectrum

$$\widetilde{F}^+(k) = \begin{cases} 2\widetilde{F}(k) = \widetilde{f}(k-k_0) + \widetilde{f}(k+k_0), & k > 0 \\ 0 & k < 0 \end{cases} \tag{18}$$

Thus, by definition, the analytic signal is a Cauchy type integral

$$F^+(s) = \frac{1}{\pi}\int_0^{\infty} \widetilde{F}(k)e^{iks}dk = \frac{1}{\pi}\int_{-\infty+i\delta}^{\infty+i\delta} F(\eta)\frac{d\eta}{s-\eta}, \tag{19}$$

regular in the lower half-plane $\mathrm{Im}\, s < 0$. For real $s$, its real part coincides with $F(s)$ whereas the imaginary part is given by Hilbert transform

$$\mathrm{Im}\, F^+(s) = \frac{1}{\pi}\,\mathrm{V.P.}\int_{-\infty}^{\infty} F(\eta)\frac{d\eta}{s-\eta} \tag{20}$$

Introduction of the analytic signal violates the causality principle: the real signal $F(ct)$ is zero before the moment of switching on the transmitter while $F^+(ct) \neq 0$ for $t < 0$. However, for a high-frequency radio pulse this discrepancy is small. So the analytic signal envelope defined as $|F^+(s)|$ is close to $f(s)$ but, contrary to the "naïve" complex signal $F_c(s)$, admits analytic continuation into complex domain, compatible with asymptotic propagation laws [22]. As an example, consider a modulated high-frequency pulse $F(s) = f(s)\cos k_0 s$ with

the envelope

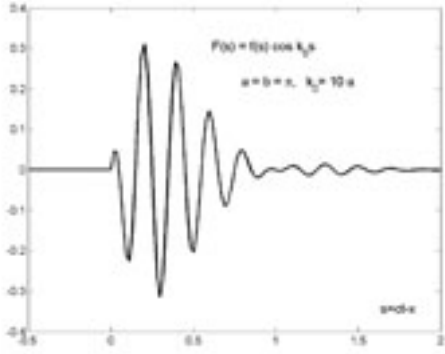 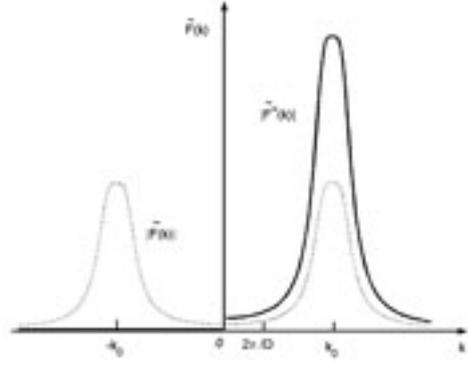

Fig. 5. Modulated radio pulse waveform (21).     Fig. 6. Analytic signal spectrum (18).

$$f(s) = \begin{cases} \sin as \exp(-bs), & s > 0 \\ 0, & s < 0 \end{cases}, \qquad \tilde{f}(k) = \frac{a}{a^2 + b^2 - k^2 - 2ibk} \qquad (21)$$

see Fig. 5. For $b \approx a$ its length is $\Lambda = cT \sim \pi/a$. The envelope spectrum $\tilde{f}(k)$ has a peak at $k = 0$ with $\Delta k \approx 2\pi/\Lambda$ and tends to zero for $|k| \to \infty$ as $O(a/k^2)$. Spectra $\tilde{F}(k)$ and $\tilde{F}^+(k)$ are shown in Fig. 6; the analytic signal envelope $|\tilde{F}^+(s)|$ is plotted in Fig. 7a,b for real and complex arguments.

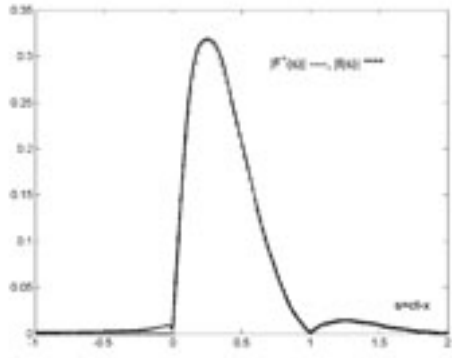 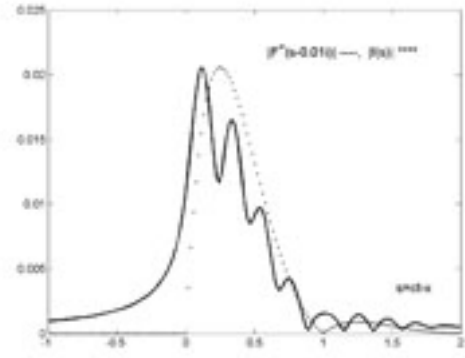

(a) (b)

Fig. 7. Analytic signal envelope of real (a) and complex (b) arguments.

For a wave packet

$$H^+(x,z,t) = \frac{1}{2\pi} \int_0^\infty \tilde{F}^+(k) u(x,z,k) e^{ik(x-ct)} dk, \qquad u(0,z,k) = A_0(z) e^{ik\Phi_0(z)} \qquad (22)$$

the "initial" condition

$$H^+(0,z,t) = \frac{A_0(z)}{2\pi} \int_0^\infty \tilde{F}^+(k) e^{ik[\Phi_0(z)-ct]} = A_0(z) F^+[ct - \Phi_0(z)] \qquad (23)$$

describes an analytic signal $F^+(ct)$ with amplitude $A_0(z)$ and initial delay $t_0(z) = \Phi_0(z)/c$.

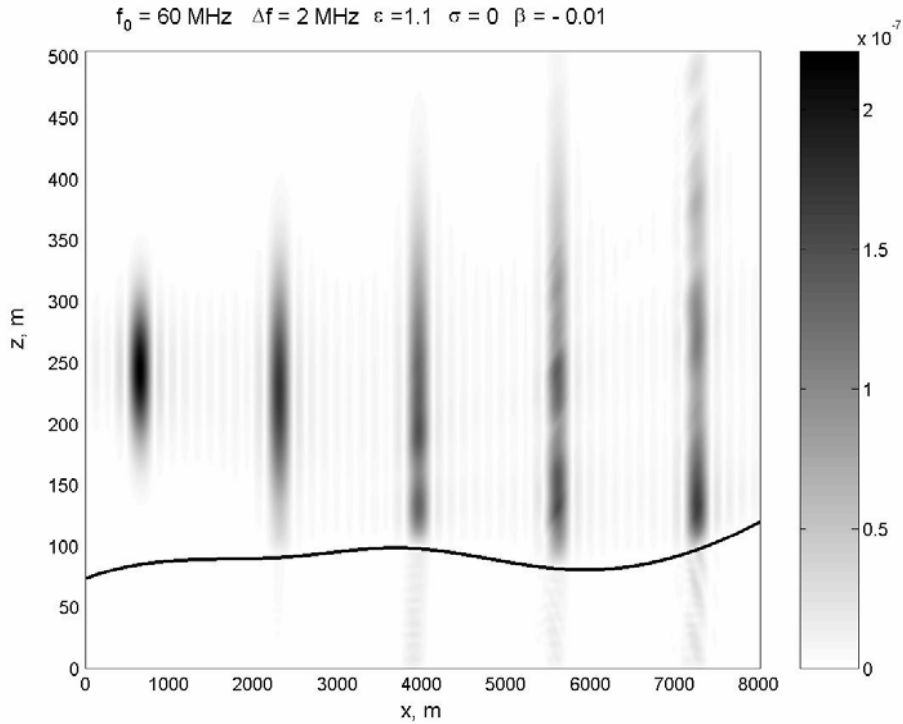

Fig. 8. Example of a Fourier-synthesized propagating EM pulse.

An example of modulated pulse propagation over smoothly rolling interface is depicted in Fig. 8. A sequence of snap-shots traces the evolution of the initial pulse envelope $f(ct)u_0(z,k_0)$ defined by (5), (21) due to the incident Gaussian beam divergence and reflection from the curved interface $z = h(x)$. It should be noted that Fourier synthesis is computationally efficient only for rather narrow-band pulses. In fact, for a good approximation of the convolution integral (22) one has to solve PE (4) for a set of wave numbers covering the spectral band $k_0 - \Delta k < k < k_0 + \Delta k$, $\Delta k \approx 2\pi/\Lambda$ with a small frequency step ($\delta k \ll \Delta k$) and, to avoid phantom solutions in the given range $\Delta x$, even more restrictive condition must be posed: $\delta k \ll 2\pi/\Delta x$. Adequate simulation methods for wide-band EM pulse propagation are discussed in the following sections.

**4. Time-domain PE and boundary conditions**

Straightforward derivation shows that if $u(x,z,k)$ is a solution of PE (4), the transient wave packet (16) $H(x,z,t) \equiv \Pi(x,z,s)$, as a function of variables $x, z, s = ct - x$, satisfies the Claerbout equation

$$2\frac{\partial^2 \Pi}{\partial x \partial s} = \frac{\partial^2 \Pi}{\partial z^2} \tag{24}$$

Equation (24), usually called "time-domain parabolic equation" (TDPE), has been obtained in [13] by formal substitution $k = i\partial/\partial s$ as well as by the reduction of the time-dependent wave equation

$$\frac{1}{c^2}\frac{\partial^2 H}{\partial t^2} = \frac{\partial^2 H}{\partial x^2} + \frac{\partial^2 H}{\partial z^2} \tag{25}$$

in a narrow vicinity of the paraxial wave front $x = ct$. Introduction of scaled variables $\xi = x/L$, $\varsigma = z/D$, $\eta = (ct-x)/\Lambda$, where $L, D$ are computational domain length and width, $\Lambda$ is spatial pulse length, yields

$$2\frac{\partial^2 \Pi}{\partial \xi \partial \eta} = \frac{\partial^2 \Pi}{\partial \varsigma^2} + \nu^2 \frac{\partial^2 \Pi}{\partial \xi^2}, \qquad \nu = D/L = \Lambda/D \ll 1 \tag{26}$$

Neglecting the small term $O(\nu^2)$ results in TDPE (24). This derivation clarifies the nature of the "time-domain parabolic equation":

1) It is a hyperbolic equation written in a traveling coordinate frame $(x,z,s)$;

2) TDPE does not describe the backward moving waves;

3) TDPE is a paraxial (narrow-angle) approximation valid in a narrow strip $D/L = O(\nu) \ll 1$;

4) TDPE describes short pulses $\Lambda/D = O(\nu) \ll 1$ whose length $\Lambda$ is comparable with the wave front deviation from the plane $x = ct$ (Fig. 9);

5) TDPE solutions are not necessary modulated high-frequency signals − they can represent short ultrawide-band pulses $f(ct)$ without carrier, e. g. a damped sinusoid (21).

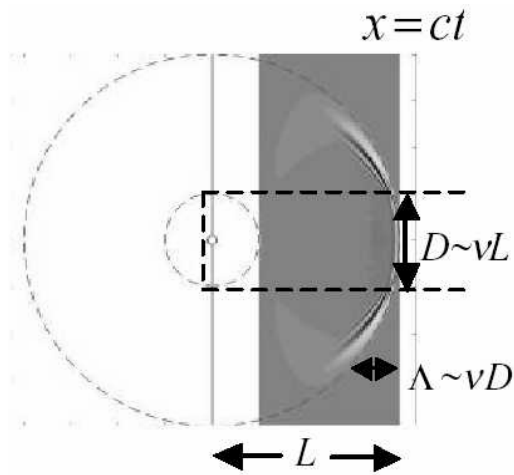

Fig. 9. Derivation of TDPE (24).

Here, a seeming contradiction may arise, as the spectral maximum of $f(s)$ can lie in the vicinity of zero frequency, not described by PE (4). As a matter of fact, at small distances from the wave front $s = O(\Lambda)$ the main part of the pulse energy is determined by the high-frequency edge of its spectrum $|k| \sim a \gg 2\pi/D$ satisfying PE applicability conditions.

To solve TDPE (24), an FD scheme of the 2nd order of accuracy has been proposed in [13]:

$$\frac{2}{\Delta x \Delta s}\left(\vec{\Pi}_{n+1,\ell+1} - \vec{\Pi}_{n,\ell+1} - \vec{\Pi}_{n+1,\ell} + \vec{\Pi}_{n,\ell}\right) = \frac{1}{4(\Delta z)^2}\nabla_z^2\left(\vec{\Pi}_{n+1,\ell+1} + \vec{\Pi}_{n,\ell+1} + \vec{\Pi}_{n+1,\ell} + \vec{\Pi}_{n,\ell}\right) \quad (27)$$

Here, $\vec{\Pi}_{n,\ell} = \{\Pi_{n,\ell}^m\}$, $x_n = n\Delta x$, $z_m = m\Delta z$, $s_\ell = \ell\Delta s$; $(\nabla_z^2 \Pi)_{n,\ell}^m = \Pi_{n,\ell}^{m+1} - 2\Pi_{n,\ell}^m + \Pi_{n,\ell}^{m-1}$. This equation is solved by zigzag marching in $(x,s)$ plane between boundary values

$$\Pi(0,z,s) = A_0(z)f[ct - \Phi_0(z)], \qquad \Pi(x,z,0) = 0 \quad (28)$$

(given source and causality condition). At each marching step $(n,m)$, a three-diagonal linear equation set arises for the unknown vector $\vec{\Pi}_{n+1,\ell+1}$. In order to complete the boundary value problem, we have to add a correct BC taking into account soil properties and to find a way of the domain truncation without creating spurious reflections. Both problems are resolved by applying Fourier transform to the frequency-domain BC (11), (13). Consider a paraxial wave packet

$$\Pi(x,z,s) = \frac{1}{2\pi}\int_{-\infty}^{\infty}\tilde{f}(k)u(x,z,k)e^{-iks}dk, \quad (29)$$

satisfying the causality condition $\Pi(x,z,s) = 0$ for $s < 0$. We rewrite the impedance BC (11) emphasizing the dependence of complex permittivity $\tilde{\varepsilon} = \varepsilon + 4\pi i\sigma/kc$ on the wave number $k = \omega/c$:

$$\frac{\partial u}{\partial z} + ik\left[\frac{\sqrt{\varepsilon-1}}{\varepsilon}\frac{\sqrt{k(k+2iq)}}{k+ir} - h'(x)u\right] = 0, \quad (30)$$

where $r = 4\pi\sigma/c\varepsilon$, $q = 2\pi\sigma/c(\varepsilon-1)$. Multiplying Eq. (30) by $\tilde{f}(k)$ and applying Fourier transform (29), we get

$$\frac{\partial \Pi}{\partial z}(x,h,s) + h'(x)\frac{\partial \Pi}{\partial s}(x,h,s) = \frac{\sqrt{\varepsilon-1}}{\varepsilon}\frac{1}{2\pi i}\int_{-\infty}^{\infty}\frac{\sqrt{k(k+2iq)}}{k+ir}\tilde{f}(k)u(x,z,k)e^{-iks}kdk. \quad (31)$$

Substituting here the inverse Fourier transform

$$\tilde{f}(k)u(x,z,k) = \int_0^{\infty}\Pi(x,z,\eta)e^{ik\eta}d\eta \quad (32)$$

we obtain, by standard calculations, the following expression for the RHS of (31):

$$\frac{\sqrt{\varepsilon-1}}{\varepsilon}\left[\frac{\partial \Pi}{\partial s}(x,h,s) - \int_0^s \frac{\partial \Pi}{\partial \eta}(x,h,\eta) N(s-\eta) d\eta\right]; \quad N(s) = e^{-rs}\left[q\int_0^s e^{(r-q)t} I_1(qt)\frac{dt}{t} + r - q\right] \quad (33)$$

Thus, we have derived a nonlocal 2D boundary condition

$$\frac{\partial \Pi}{\partial z}(x,h,s) + h'(x)\frac{\partial \Pi}{\partial s}(x,h,s) = \frac{\sqrt{\varepsilon-1}}{\varepsilon}\left[\frac{\partial \Pi}{\partial s}(x,h,s) - \int_0^s \frac{\partial \Pi}{\partial \eta}(x,h,\eta) N(s-\eta) d\eta\right] \quad (34)$$

being an exact time-domain counterpart of the impedance BC (11). Its nonlocality is a consequence of interaction between two waves propagating along the earth surface with different phase velocities. The integral term kernel $N(s-\eta)$ can be easily calculated for different $\varepsilon$ and $\sigma$. For $\varepsilon > 3.15$, function $N(s)$ monotonously tends to zero with increasing $s$ – see Fig. 10.

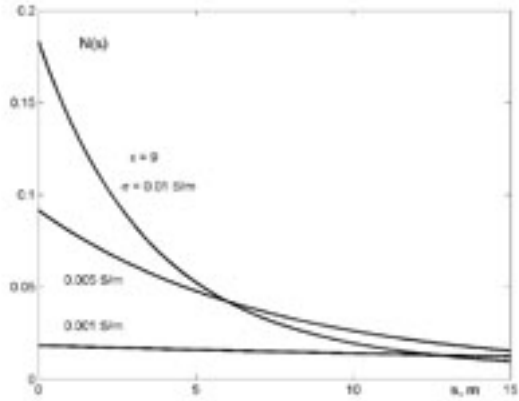

Fig. 10. Kernel of nonlocal impedance BC (34).

It is interesting to note that in both limiting cases: perfectly conducting boundary ($\sigma = \infty$) and zero soil conductivity ($\sigma = 0$) Eq. (34) reduces to a local BC. In the former case $N(s)$ tends to a delta function, for $N(0) = r - q \underset{\sigma \to \infty}{\to} \infty$, $\int_0^\infty N(s) ds = 1$. The integral in (34) limits to $\frac{\partial \Pi}{\partial s}(x,h,s)$, the RHS vanishes, and we get a Neumann BC. In the opposite case ($\sigma \to 0$) the spatial scale of $N(s)$ (length of the pulse "dispersion tail") is growing but its absolute value is tending to zero, so only a local term $\frac{\sqrt{\varepsilon-1}}{\varepsilon}\frac{\partial \Pi}{\partial s}(x,z,s)$ remains.

In a similar way, the time-domain generalization of the transparency BC (13) is derived which grants the absence of spurious reflections a from the artificial computational boundary $z = z_{max}$. Applying to Eq. (13) the Fourier transform (29), (32) and denoting $k = ip$, we get

$$\frac{\partial \Pi}{\partial z}(x, z_{max}, s) = \frac{i}{\pi\sqrt{2\pi}} \int_{-i\infty}^{i\infty} e^{ps} \sqrt{p}\, dp \int_0^x \frac{d\xi}{\sqrt{x-\xi}} \int_0^\infty e^{-p\eta} \frac{\partial \Pi}{\partial \xi}(\xi, z_{max}, \eta)\, d\eta =$$
$$= \frac{i}{\pi\sqrt{2\pi}} \int_0^x \frac{d\xi}{\sqrt{x-\xi}} \int_0^\infty \frac{\partial^2 \Pi}{\partial \xi \partial \eta}(\xi, z_{max}, s) \int_{-i\infty+0}^{i\infty+0} e^{p(s-\eta)} \frac{dp}{\sqrt{p}} \quad (35)$$

Finally, evaluating the inner integral we obtain an elegant 2D boundary condition

$$\frac{\partial \Pi}{\partial z}(x, z_{max}, s) = -\frac{\sqrt{2}}{\pi} \int_0^x \int_0^s \frac{\partial^2 \Pi}{\partial \xi \partial \eta}(\xi, z_{max}, \eta) \frac{d\xi d\eta}{\sqrt{(x-\xi)(s-\eta)}} \quad (36)$$

symmetric with respect to the variables $x, z$, which could be expected from the symmetry of the TDPE (24).

A simulated example of ultrawide-band EM pulse propagation over a nonuniform earth surface with soil parameters $\varepsilon = 10$, $\sigma = 9 \cdot 10^7 s^{-1}$ (0.01 $S/м$) is depicted in Fig. 11a. Evolution of the spatial amplitude distribution for a pulsed signal generated by a Gaussian source: $A_0(z) = \exp[-(z-z_0)^2/w_0^2]$ with a skew curved wave front: $\Phi_0(z) = (z-z_0)^2/2\rho_0 + \beta(z-z_0)$ is shown in a grey color scale. The initial pulse waveform $f(ct)$ is a damped sinusoid (21) with $a = b$ and spatial length $\Lambda = \pi/a = 30\,m$.

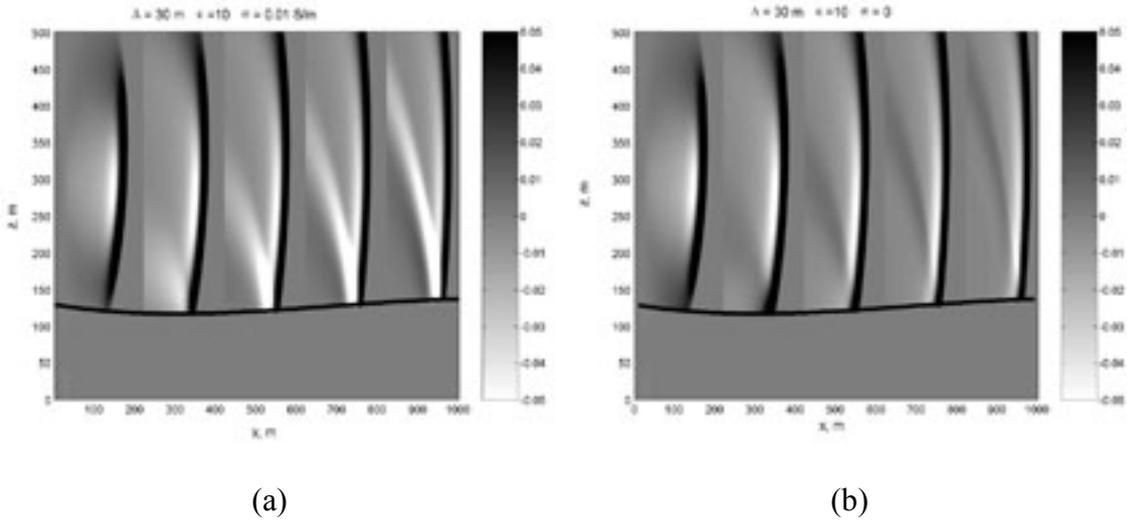

(a)             (b)

Fig. 11. Propagation of ultrawide-band pulse (21) over nonuniform earth surface.
Initial Gaussian beam parameters: $z_0 = 300\,m$, $w_0 = 80\,m$, $\rho_0 = 300\,m$, $\beta = -0.1$.
Soil conductivity: $\sigma = 0.01\,S/m$ (a), $\sigma = 0$ (b).

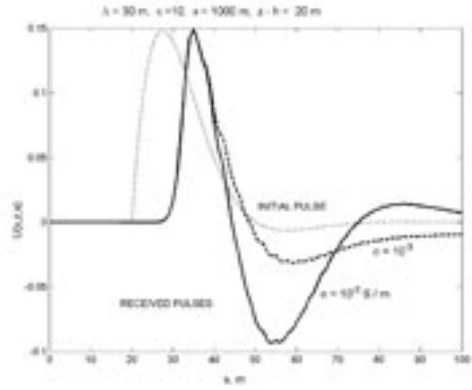

Fig. 12. Received pulse waveform depending on soil conductivity:
$\sigma = 0.01\ S/m$ (solid line), $\sigma = 0.001\ S/m$ (dashed line); initial pulse (dots).

The snapshots clearly show the reflected pulse generation at the earth surface. The transparency BC (36) imposed at the height $z_{max} = 500\,m$ assures unimpeded radiation exit from the computational domain. Finite soil conductivity causes signal dispersion appearing in a certain delay of the reflected pulse. It is obvious from the comparison with a similar plot calculated for a model non-conducting soil $\sigma = 0$ - Fig. 11b. A quantitative estimation of the effect can be made by means of Fig. 12 revealing a considerable dependence of the pulse waveform on the soil conductivity. This effect caused by the pulse penetration into the ground can be used for ecological monitoring (water pollution, earthquake precursors, etc.).

## 5. Hybrid TDPE and short high-frequency pulse propagation

An important practical issue is overland propagation of short EM pulses with high-frequency carrier. Basically, having absolute stability, TDPE (24) is capable to describe wide-band radio pulse propagation. However, the computational expense is drastically growing with increasing carrier frequency. If $\omega_0$ considerably exceeds the spectral band of the signal it is useful to factor out the carrier:

$$\Pi(x,z,s) = U(x,z,s)\exp[i(k_0 x - \omega_0 t)] \tag{37}$$

and to consider the transient signal envelope $U(x,z,s)$ satisfying a hybrid equation

$$2\left(ik_0 - \frac{\partial}{\partial s}\right)\frac{\partial U}{\partial x} + \frac{\partial^2 U}{\partial z^2} = 0, \tag{38}$$

combining the features of standard Leontovich PE (4) with Claerbout TDPE (24). Unfortunately, despite a relatively slow variations of $U(x,z,s)$ in space-time, straightforward numerical solution of this hybrid TDPE (HPE) entails considerable difficulties, as the large coefficient $k_0$ by $\partial U/\partial x$ demands a dense computational grid. On the other hand, the presence of a big parameter allows us to construct an asymptotic solution of HPE (38) radically reducing the computational burden.

In order to find a proper asymptotic Ansatz, consider monochromatic PE (4) at the carrier frequency $\omega_0 = k_0 c$ with the initial condition $u(0,z) = A_0(z) e^{ik_0 \Phi_0(z)}$. We admit complex values of the eikonal $\Phi_0(z)$ to describe relatively narrow wave packets, like a Gaussian beam (6). For $k_0 \to \infty$, we obtain an asymptotic solution

$$u(x,z,k_0) = A(x,z) e^{ik_0 \Phi(x,z)} \tag{39}$$

where $\Phi(x,z)$ satisfies a "parabolic" eikonal equation

$$\Phi_x + \frac{1}{2}\Phi_z^2 = 0 \tag{40}$$

while slowly varying amplitude $A(x,z)$ is governed by the paraxial transport equation

$$A_x + \Phi_z A_z + \frac{1}{2}\Phi_{zz} A = 0 \tag{41}$$

Eqs. (40)-(41), being an approximate form of the well-known laws of geometric optics (GO) [4], can be easily solved by the method of characteristics. Consider a particular solution of the eikonal equation (40) corresponding to a bundle of rays spreading from a central point $x = 0$, $z = z_0$:

$$\Phi(x,z) = \frac{(z-z_0)^2}{2x} \approx \sqrt{x^2 + (z-z_0)^2} - x \tag{42}$$

At a characteristic line $z = z_0 + \gamma x$ we have $\Phi(x, \gamma x) = \frac{1}{2}\gamma^2 x$. The envelope of the family (42) solves the boundary value problem with an arbitrary initial condition $\Phi(0,z) = \Phi_0(z)$. Define

$$\Phi[x, z_0 + \gamma(z_0)x] = \Phi_0(z_0) + \frac{1}{2}\gamma^2(z_0) x \tag{43}$$

By differentiating Eq. (43) with respect to $x$ and $z_0$ we obtain

$$\Phi_x + \gamma \Phi_z = \frac{1}{2}\gamma^2, \qquad (1+\gamma' x)\Phi_z = \Phi_0' + \gamma\gamma' x \tag{44}$$

Function (43) will satisfy the eikonal equation (40) if the ray direction is matched with the local wave front tilt: $\gamma(z_0) = \Phi_0'(z_0)$. Having constructed the eikonal $\Phi(x,z)$ we reduce the transport equation (41) to an ODE

$$\frac{d}{dx}A(x,z_0+\gamma x)+\frac{1}{2}\frac{\gamma'(z_0)}{1+\gamma'(z_0)x}A(x,z_0+\gamma x)=0 \qquad (45)$$

with an evident integral

$$A[x,z_0+\gamma(z_0)x]=\frac{A_0(z_0)}{\sqrt{1+\gamma'(z_0)x}} \qquad (46)$$

In a similar way, an asymptotic solution of the modified Claerbout equation (38) can be found. Substituting the Ansatz $U(x,z,s) = I(x,z,s)\exp[ik_0\Phi(x,z)]$ into Eq. (38) we obtain

$$-k_0^2\left(\Phi_x+\frac{1}{2}\Phi_z^2+ik_0\right)+ik_0\left(I_x-\Phi_x I_s+\Phi_z I_z+\frac{1}{2}\Phi_{zz}I\right)+\frac{1}{2}I_{zz}-I_{xs}=0 \qquad (47)$$

The leading term $O(k_0^2)$ disappears in virtue of the eikonal equation (40). Thus, to the accuracy $O(k_0^{-1})$, a space-time transport equation arises for the slowly varying amplitude $I(x,z,s)$:

$$I_x-\Phi_x I_s+\Phi_z I_z+\frac{1}{2}\Phi_{zz}I=0 \qquad (48)$$

As the coefficients of Eq. (48) do not depend on *s*, it has a solution of the following form

$$I(x,z,s)=A(x,z)g[s-\Psi(x,z)] \qquad (49)$$

Here, $A(x,z)$ is a solution of the stationary transport equation (41) while $g(s)$ is an arbitrary function of $s=ct-x$, and $\Psi(x,z)$ satisfies a linear PDE

$$\Psi_x+\Phi_z\Psi_z=-\Phi_x \qquad (50)$$

Solving Eq. (50) by characteristics one easily gets

$$\Psi[x,z_0+\gamma(z_0)x]=\Phi[x,z_0+\gamma(z_0)x]+\theta(z_0) \qquad (51)$$

where $\theta(z_0)$ is an arbitrary function. So, the solution of the HPE (38) has asymptotic representation

$$U(x,z,s)\underset{k_0\to\infty}{\sim}A(x,z)g[s-\Psi(x,z)]\exp[ik_0\Phi(x,z)]=g[s-\Phi(x,z)-\theta(z_0)]u(x,z,k_0) \qquad (52)$$

Here, $u(x,z,k_0)$ is a solution of the standard Leontovich PE (4), $A(x,z)$ and $\Phi(x,z)$ are its amplitude and eikonal, respectively; $g(s)$ and $\theta(z_0)$ are arbitrary functions, and $z_0(x,z)$ is to be found from the transcendental equation $z_0+\Phi_0'(z_0)x=z$. Asymptotic solution (52) is a paraxial version of the space-time GO [29], the rays and wave fronts being defined numerically via parabolic equation. Functions $A(x,z)$ and $\Phi(x,z)$ are generally complex-valued, so distinction between wave amplitude and complex "phase" is made solely on the basis of their different dependence on frequency. In particular, complex eikonal $\Phi(x,z)$, defined as

$$\Phi(x,z) = -i\frac{\partial}{\partial k_0}\log u(x,z,k_0) \approx \frac{u(k_1) - u(k_2)}{i(k_1 - k_2)u\left(\frac{k_1 + k_2}{2}\right)} \tag{53}$$

is calculated from PE numerical solutions at two close frequencies $\omega_{1,2} = k_{1,2}c$.

Physically, complex eikonal in Eq. (52) appears due to diffraction effects described by PE (4). An important consequence is the absence of singularities in the constructed asymptotic solution, as the "parabolic" rays do not produce caustics in the real space. Another effect caused by diffraction – pulse envelope distortion also is taken into account via complex values of the signal delay $\psi(x,z)/c$. Physical meaning of complex $s = ct - x$ is provided by the theory of analytic signal [22].

Arbitrary functions in Eq. (52) are determined by the initial conditions. In the simplest case the constructed transient (37) has the form

$$H(x,z,t) \equiv \Pi(x,z,ct-x) \approx A(x,z)F^+[ct - x - \Phi(x,z)] \tag{54}$$

where $A(x,z)$ and $\Phi(x,z)$ are complex amplitude and eikonal evolved from the initial $A_0(z)$, $\Phi_0(z)$ given by Eq. (23) and $F^+(ct) \approx f(ct)\exp(-i\omega_0 t)$ is the analytic signal (19) corresponding to the real signal $f(ct)$. Physical meaning has the real part of the complex solution (54) or, from the practical point of view, its normalized envelope $|H(x,z,s)|/\sqrt{2}$.

In virtue of the superposition principle, a more general asymptotic solution can be constructed as a number of terms (54). That is a direct analogy with ordinary GO where the incident and reflected waves correspond to different ray families. An important practical example is radar pulse propagation over the earth surface when the direct and reflected from the ground pulses can be distinguished and used for target location [23].

Consider first a model example: a short pulse with carrier frequency $f_0 = 200$ MHz and damped sinusoidal envelope (21), propagating over a slowly rolling boundary $z = h(x)$. Initial pulse parameters are: $a = b$, $\Lambda = \pi/a \approx 9\ m$; $w_0 = 15\ m$, $\beta = -0.01$, $\rho_0 = 200\ m$. Stationary field distribution calculated by numerical integration of PE (4) at $k_0 = 2\pi f_0/c$ produces a regular interference pattern – Fig. 13. It can be represented as a superposition of the incident Gaussian beam (6) $u_i(x,z,k_0)$ with the reflected wave $u_r(x,z,k_0) \equiv u - u_i$ determined by the terrain $z = h(x)$ and the impedance BC (11). Functions $A_i(x,z), \Phi_i(x,z)$ are given by the asymptotic solution (39)-(41), and eikonal $\Phi_r(x,z)$ is reconstructed from the spatial phase distribution of the reflected wave – see Fig. 14. In accordance with such a monochromatic framework, a two-term asymptotic formula arises for the pulsed transient:

$$H(x,z,t) \sim A_i(x,z) F^+[ct - x - \Phi_i(x,z)] + A_r(x,z) F^+[ct - x - \Phi_r(x,z)]$$

(55)

$$= u_i(x,z) e^{-ik\Phi_i(x,z)} F^+[ct - x - \Phi_i(x,z)] + u_r e^{-ik\Phi_r(x,z)}(x,z) F^+[ct - x - \Phi_r(x,z)]$$

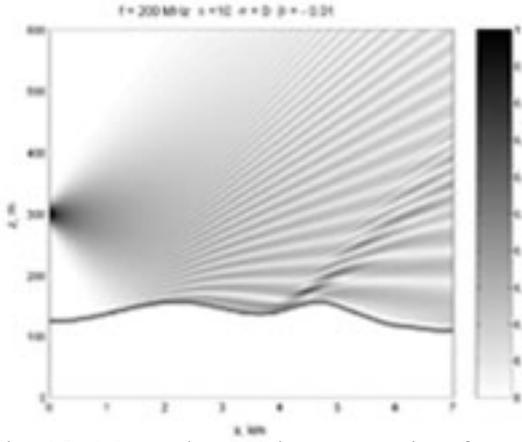 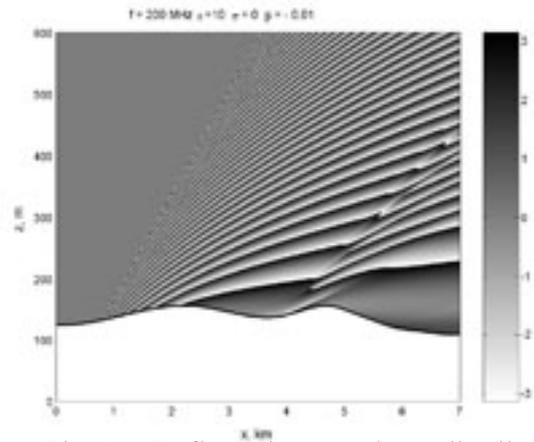

Fig. 13. Monochromatic attenuation function.    Fig. 14. Reflected wave phase distribution.

Note that to find the amplitude and complex delay of the incident and reflected signals we need just to solve the standard PE in frequency domain at two close frequencies $\omega_{1,2} \approx \omega_0$ - see Eq. (53). In time domain, calculation reduces to the evaluation of an analytic function $F^+(s)$ for the given argument values of interest. That radically diminishes the required computational resources compared with direct numerical integration of the TDPE (24).

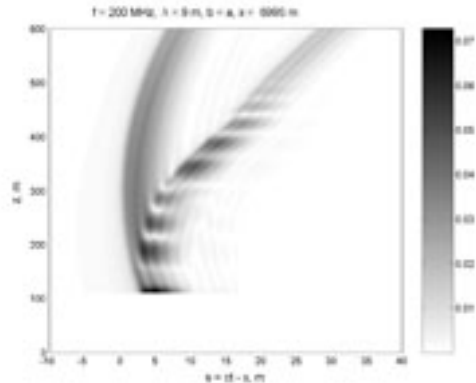

Fig. 15. Received pulse envelope as function of receiver altitude.

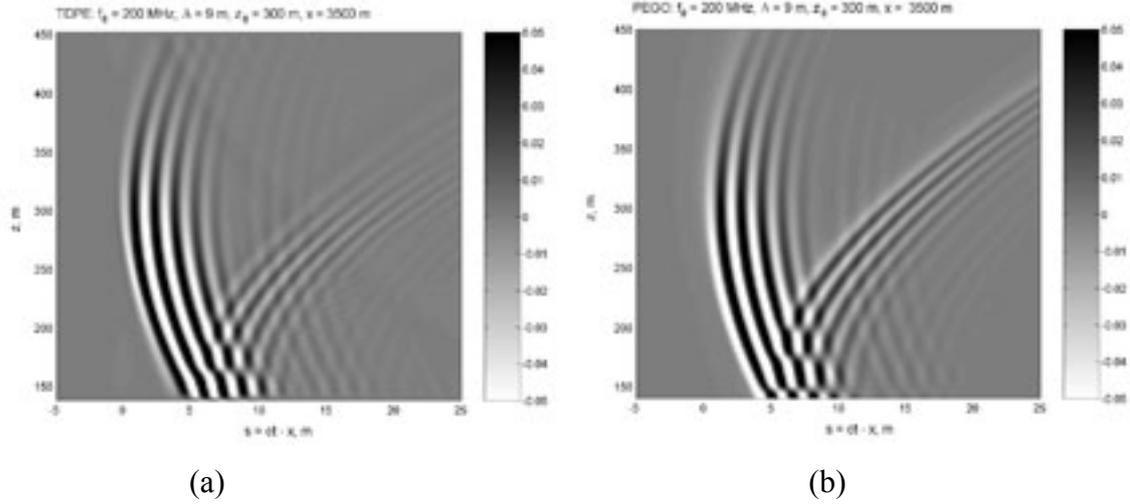

(a)  (b)

Fig. 16. Comparison between numerical solution of TDPE (24) (a) and asymptotic solution (55) (b).

The envelope of the received analytic signal (55), as a function of $s = ct - x$ and the receiver height $z$, for $x = 7\,km$ is shown in Fig. 15. One can see profound interference minima near the earth surface and a good separation of the direct and reflected pulse to heights above 400 $m$. Figures 16a,b compare the asymptotic solution with direct numerical integration of the Claerbout TDPE (24). Qualitatively, they are almost identical. Some hardly seen discrepancy is due to a limited accuracy of the asymptotic solution and FD scheme (27). This comparison demonstrates the efficiency of the developed approach. Substantial acceleration of the numerical procedure (by around 200 times in this example) makes a good reason to use it in realistic conditions.

As an example, we simulate an experimental situation [23]: radar pulse propagation between two aircrafts flying by parallel routs over an irregular terrain. The experiment [23] was aimed at simultaneous determination of the target range and altitude from the measured return times of the direct radar pulse and the echo signal from the earth surface. Our goal is to develop an efficient method of EM field calculations under conditions of multipath and signal distortion. At such large ranges ($X = 100\,km$) the Earth sphericity must be taken into account. For this purpose, a parabolic hump $x(X - x)/2R^*_{earth}$ has been added to the real terrain profile plotted in .9 of [23]. Atmospheric refraction has been considered by using the equivalent Earth radius $R^*_{earth} = 4/3\,R_{earth}$ [12]. Global field strength distribution produced by the incident carrier wave at $f_0 = 141\,\mathrm{MHz}$ is depicted in Fig. 17a. Despite evident multipath character of the reflected wave (Fig. 17b) its eikonal $\Phi_r(x,z)$ has a rather regular structure. Therefore, our PE based version of complex GO can be applied to simulate the averaged parameters of the received radar pulse (the actually observed signal is a stochastic quantity with normal distribution [23]). Its envelope, as a function of the distance from the paraxial wave front $s = ct - x$ and the receiver height $z$, is depicted in Fig. 18.

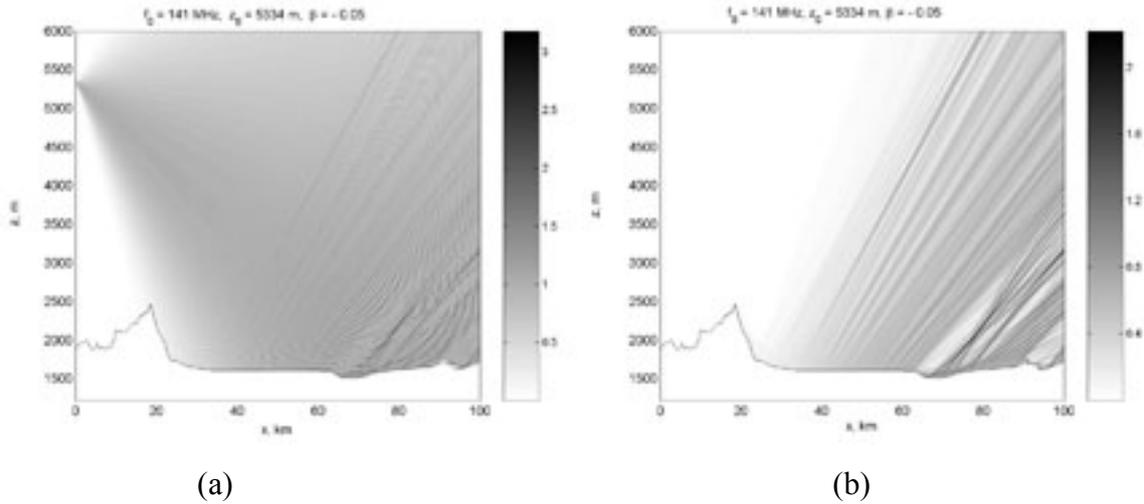

(a)                                   (b)

Fig. 17. Simulation of experiment [23]: global field strength (attenuation function) at $f_0 = 141$ $MHz$ $(a)$, reflected wave (b).

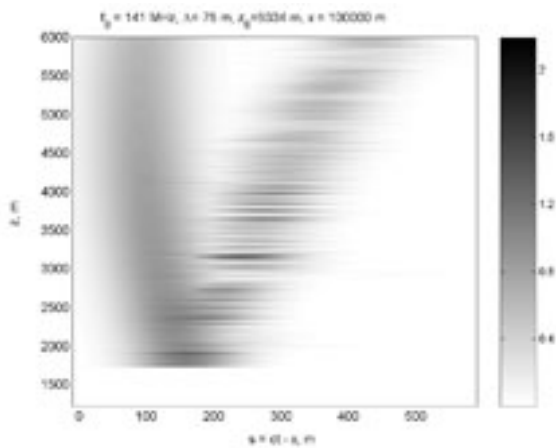   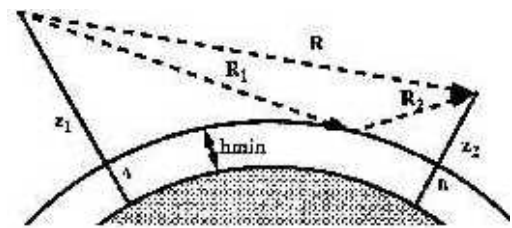

Fig. 18. Gaussian pulse envelope as function     Fig. 19. Target height from reflected
of relative delay and receiver altitude.                pulse delay [23].

A Gaussian pulse waveform is chosen with the parameters corresponding to the experimental situation [23]: $z_0 = 5.3$ $km$, $f_0 = 141$ MHz, $\Lambda \approx 75$ $m$. The direct and reflected pulses are distinctly separated for $z > 4.5 km$ which allows one to reliably solve the triangle $R, R_1, R_2$ for target altitude determination - see Fig. 19 borrowed from [23]. Variability and statistics of the simulated reflected pulse resemble the experimental plots presented in [23], and the calculated received pulsed signal envelope for a fixed receiver height $z = 5.3\ km$ (Fig. 20) agrees well both with the experimental data and the results of thorough statistical modeling [23]. The quantitative discrepancy in the reflected pulse amplitude does not exceed the inherent uncertainty due to the errors in terrain description.

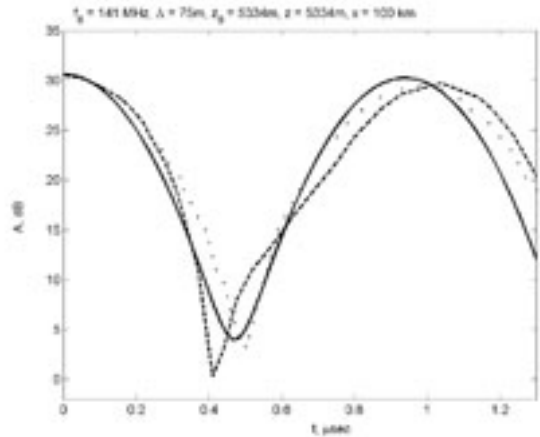

Fig. 20. Received pulse envelope: asymptotic HPE solution (solid line), experiment [23] (dash), stochastic model calculation [23] (dots).

**Acknowledgements**

This work was supported in part by a joint RFBR-DFG grant No 01-02-04003. The authors are grateful to Friedrich Landstorfer and Ningyan Zhu who initiated this research.